\documentclass[conference]{IEEEtran}      

\IEEEoverridecommandlockouts                              

\usepackage{amsfonts}
\usepackage{amsmath,amsthm,amssymb,enumerate}
\usepackage{graphicx}
\usepackage{caption}
\usepackage{epstopdf}
\usepackage{algorithm}
\usepackage{algorithmicx}

\usepackage{verbatim}
\usepackage{algpseudocode}
\usepackage{mathtools}
\usepackage{flushend}
\usepackage{setspace}

\usepackage{etoolbox}

\usepackage{xcolor}

\usepackage[utf8]{inputenc}
\usepackage{url}
\usepackage{graphics}
\usepackage{graphicx}
\usepackage{amsfonts}
\usepackage{amsmath,amsthm,amssymb,enumerate}
\usepackage{mathrsfs}
\usepackage{cite}
\usepackage{enumerate}
\usepackage{comment}
\usepackage{subfig}
\usepackage{multirow}
\usepackage[english]{babel}
\usepackage{ifpdf}
\usepackage{color}
\usepackage{nicefrac}
\usepackage{breqn}
\usepackage{pgfplots}

\newcommand{\prox}{\mathop{\rm prox}\nolimits}

\newcommand{\sign}{\mathop{\rm sign}\nolimits}

\newcommand{\minimize}{\operatorname{minimize}}  
\renewcommand{\Re}{{\rm{I\!R}} }

\DeclareMathOperator*{\argmin}{\arg\!\min}

\definecolor{red}{rgb}{1,0,0}
\definecolor{yellow}{rgb}{1,1,0}
\definecolor{blue}{rgb}{0,0,1}

\newcommand{\hide}[1]{}

\newcommand{\R}{\mathbb{R}}

\newcommand{\T}{\top}

\newtheorem{rem}{Remark}

\DeclareMathOperator{\supp}{supp}

\title{\LARGE \bf Accelerated reconstruction of a compressively sampled data stream}

\author{\IEEEauthorblockN{Pantelis Sopasakis\IEEEauthorrefmark{1},
Nikolaos Freris\IEEEauthorrefmark{2}, and Panagiotis Patrinos\IEEEauthorrefmark{3}
 }

\IEEEauthorblockA{\IEEEauthorrefmark{1}IMT Institute for Advanced Studies Lucca,
Piazza S. Ponziano 6, 55100 Lucca, Italy
}
\IEEEauthorblockA{\IEEEauthorrefmark{2}Division of Engineering, New York University Abu Dhabi, P.O. Box 129188, Abu Dhabi, United Arab Emirates}
\IEEEauthorblockA{\IEEEauthorrefmark{3}STADIUS Center for Dynamical Systems, Signal Processing and Data Analytics,\\
	KU Leuven, Department of Electrical Engineering (ESAT), Kasteelpark Arenberg 10, 3001 Leuven, Belgium}

\thanks{The work of the first author was supported by the EU-funded H2020 research project DISIRE, grant agreement No. 636834.}
}

\begin{document}
\maketitle

\begin{abstract}
	The traditional compressed sensing approach is naturally offline, in that it amounts to sparsely sampling and reconstructing a given dataset. Recently, an online algorithm for performing compressed sensing on streaming data was proposed~\cite{RCS}: the scheme uses recursive sampling of the input stream and recursive decompression to accurately estimate stream entries from the acquired noisy measurements. 
	
	In this paper, we develop a novel Newton-type forward-backward proximal method to recursively solve the regularized Least-Squares problem (LASSO) online. We establish global convergence of our method as well as a local quadratic convergence rate. Our simulations show a substantial speed-up over the state of the art which may render the proposed method suitable for applications with stringent real-time constraints.	
\end{abstract}

\section{Introduction}\label{sec:Introduction}

In signal processing, continuous signals are typically sampled at discrete instances in order to store, process, and share. In doing so, signals of interest are assumed sparsely representable in an appropriately selected orthonormal basis, i.e., they can be reconstructed by storing few non-zero coefficients in the given basis. For instance, the Fourier basis is used for bandlimited signals -- with prominent applications in communications -- while wavelet bases are suitable for representing piecewise smooth signals, such as bitmap images. Traditionally, the celebrated Nyquist-Shannon sampling theorem suggests a sampling rate that is at least twice the signal bandwidth; however this rate may be  unnecessarily high compared to the signal's \emph{innovation}, i.e., the minimum number of coefficients sufficient to accurately represent it in an appropriately selected basis.

\emph{Compressed Sensing} (CS)~\cite{Donoho2006,Candes2006b} is a relatively new sampling paradigm that was introduced for sampling signals based entirely on their innovation, and has become ever since a major field of research in signal processing and information theory. The major contribution of this framework is a lower sampling rate compared to the classical sampling theory for signals that have sparse representation in some fixed basis~\cite{Candes2008}, with notable applications in imaging~\cite{CS_applications}, such as MRI. 

\subsection{Compressed Sensing}
For a vector $x \in \R^n$ we define the $\ell_0$ pseudo-norm as the cardinality of its support $||x||_0 := |\mbox{supp}(x)|$, where the support is the set of non-zero entries $\supp (x) := \{i: x_i \neq 0\}$. A vector $x$ is $s$\emph{-sparse} if and only if it has at most $s$ non-zero entries, $||x||_0\le s$. 

\textbf{Compression:} CS performs linear sampling $y=Ax$, where $A\in\R^{m \times n}$. Compression is performed by obtaining $m \ll n$ measurements. The main result states that \emph{any} $s$-sparse vector $x\in \R^n$, can be sampled by a \emph{universal}\footnote{The matrix A is universal in the sense it may be used to sample any vector with no more than $s$ non-zeroes, regardless of their positions.} matrix $A$ with $m = \Theta(s \log(\frac{n}{s}))$, and may then be reconstructed perfectly solely from (noiseless) measurements $y$~\cite{RCS}. For this purpose, matrix $A$ needs to satisfy a certain \emph{Restricted Isometry Property}: there exists $\delta_s \in [0,1)$ sufficiently small so that: 
\begin{align*}
(1-\delta_s)||x||_2^2 \le ||Ax||_2^2 \le (1+\delta_s)||x||_2^2
\end{align*}
holds for all $s$-sparse vectors $x$. For many practical applications, taking $m=4s$ works well. The success of CS lies in that the sampling matrix can be generated very efficiently as a random matrix, e.g., having i.i.d. Gaussian entries $\mathcal{N}\left(0,1/m\right)$ .
Note that in practice measurements are typically noisy, $y=Ax+w$. 

\textbf{Decompression:} In order to retrieve the original vector $x$ from noisy measurements $y=Ax+w$, one needs to solve the $\ell_1$-regularized least squares problem: 
\begin{align}
\minimize&\  \tfrac{1}{2}\Vert Ax-y \Vert_2^2 + \lambda\Vert x \Vert_1,
\label{eq:lasso}
\end{align}
where $\lambda$ is the regularization parameter that controls the trade-off between sparsity and reconstruction error. This is best known as Least Absolute Selection and Shrinkage Operator (LASSO) in the statistics literature \cite{Tibshirani1996}. There are several results analyzing the reconstruction accuracy of LASSO; for example~\cite{candes2009} states that if $w \sim \mathcal{N}(0,\sigma^2I)$,  $\min_{i\in \supp(x)} |x_i|  > 8 \sigma \sqrt{2 \log n}$, and we choose $\lambda = 4 \sigma\sqrt{2 \log n}$ then a solution $x^{\star}$ to~\eqref{eq:lasso}  has the same support as $x$, and its non-zero entries the same sign with their corresponding ones of $x$, 
with high probability. Additionally, the $\ell_2$ reconstruction error is proportional to  the standard deviation of the noise $\sigma$~\cite{Candes2008}. 

\begin{rem}[Algorithms for LASSO] 
	LASSO can easily be recast as a quadratic program which can be handled by interior point methods~\cite{kim2007interior}.
	Additionally, iterative algorithms have been developed specifically for LASSO; all these are inspired by proximal methods  \cite{Parikh:2013} for non-smooth convex optimization: FISTA \cite{Beck2009} and SpaRSA \cite{Wright:2009kw} are accelerated proximal gradient methods \cite{Parikh:2013}, SALSA \cite{Afonso2010} is an application of the alternative direction method of multipliers (ADMM). These methods are \emph{first-order} methods, in essence generalizations of the gradient method and feature \emph{sublinear convergence}.    
In the current paper, we devise a proximal Newton-type method 
	with substantial speedup exploiting the fact that its convergence rate is locally \emph{quadratic} (i.e., goes to zero roughly like $e^{-ct^2}$ at the vicinity of an optimal solution).
\end{rem}


\section{Recursive Compressed Sensing}\label{sec:RCS}

\subsection{Problem statement}
The traditional CS framework is naturally \emph{offline} and requires compressing and decompressing an entire given dataset at one shot. The Recursive Compressed Sensing~\cite{RCS, allerton} was developed as a new method for performing CS on an infinite data stream. The method consists in successively sampling the data stream via applying traditional CS to sliding overlapping windows in a recursive manner. 
Consider an infinite sequence, $\{x_i\}_{i = 0,1,\dots}$. Let
$$S:= \lim\sup_{j\to\infty} j^{-1}|\supp(\{x_i\}_{i=0,\hdots,j-1})| \in [0,1]$$
be the \emph{average} sparsity. 
We define successive windows of length $n$:
 \begin{align}
 x^{(i)} := \begin{bmatrix}
 x_{i} &x_{i+1} & \dots &x_{i+n-1}
 \end{bmatrix}^{\T}
 \label{eq:wstructure}
 \end{align}
and take $s = Sn$ the average sparsity parameter for a window of length $n$.
The sampling matrix $A\in \R^{m\times n}$ (where we may take $m = \Theta(s \log(\frac{n}{s}))$) is only generated once. Each window is compressively sampled given a matrix $A^{(i)}$:
$$y^{(i)} = A^{(i)}x^{(i)} + w^{(i)}.$$
The sampling matrix $A^{(i)}$ is recursively computed by $A^{(0)} = A$ and 
$A^{(i)} = A^{(i-1)}P,$
where $P$ is a permutation matrix which when right-multiplying a matrix cyclically rotates its columns to the left. This gives rise to an efficient recursive sampling mechanism for the input stream, where measurements for each window are taken via a rank-1 update~\cite{RCS}. 

For \emph{decompression}, we need to solve the LASSO for each separate window:
\begin{align}
\minimize&\  \tfrac{1}{2}\Vert A^{(i)}x^{(i)}-y^{(i)} \Vert_2^2 + \lambda\Vert x^{(i)} \Vert_1.
\end{align}

Note that the windows overlap, hence for a given stream entry, multiple estimates (one per each window that contains it) may be obtained. These estimates are then combined to boost estimation accuracy using non-linear support detection, least-squares debiasing and averaging with provable performance amelioration~\cite{RCS}. 

The overlap in sampling can further be exploited to speed-up the stream reconstruction: we use the estimate from a previously decompressed window to warm-start the numerical solver for LASSO in the next one. This simple idea provides a mechanism for efficient \emph{recursive estimation}. 

Formally, let  
$x^{(i)} =  [ \hat{x}_{0}^{(i)}\ \dots\ \hat{x}_{n-1}^{(i)} ]^{\top}$
denote the optimal estimate obtained by LASSO in the $i$-th window, where we use $\hat{x}_{j}^{(i)}$ to denote the $j$-th entry of the $i-$th window (which according to our definitions corresponds to the $(i+j)$-th stream entry). In order to solve LASSO to obtain $\hat{x}^{(i+1)}$ we may use: 
\begin{align*}
	\hat{x}_{[0]}^{(i)} 
	&= \begin{bmatrix}
		\hat{x}_1^{(i-1)} & \hat{x}_2^{(i-1)} & \dots & \hat{x}_{n-1}^{(i-1)} & 0 
	\end{bmatrix}^\T,
\end{align*}
as the starting point in the iterative optimization solver for LASSO in the $(i+1)-$th window, where $\hat{x}^{(i-1)}_j$, for $j=1,\dots,n-1$, is the portion of the optimal solution based on the previous window. 
The last entry $\hat{x}^{(i)}_{n-1}$ is set to 0, since we are considering sparse streams with most entries being zero. 

\subsection{Our contribution}
We devise a new numerical algorithm for solving LASSO based on the recently developed idea of proximal envelopes~\cite{patrinos2013proximal,patrinosFBN}. The new method demonstrates favorable convergence properties when compared to first order methods (FISTA, SpaRSA, SALSA, L1LS~\cite{kim2007interior});  in particular it has local quadratic convergence. Furthermore, this scheme is very efficient as each iteration boils down to solving a linear system of low dimension. Using this solver in RCS along with warm-starting leads to substantive acceleration of stream decompression. We verify this with a rich experimental setup. Unfortunately, due to length constraints, we do not provide a complete convergence analysis, but rather some hints on global convergence and local quadratic rate. 

\section{Algorithm}\label{sec:Algorithm}
\subsection{Forward Backward Newton Algorithm}
Let $f(x)=\tfrac{1}{2}\|Ax-y\|_2^2$, $g(x)=\lambda\|x\|_1$. Then $x^\star$
is optimal for
\begin{equation}\label{eq:str_prb}
\minimize\ \varphi(x)=f(x)+g(x)
\end{equation}
if and only if it satisfies
\begin{equation}\label{eq:gen_optCond}
-\nabla f(x^\star)\in\partial g(x^\star)
\end{equation}
where $\partial g(x)$ is the subdifferential of $g$ at $x$, defined by:
\begin{align}
 \partial g(x){=}\left\{v\in \Re^n\mid g(w)\geq g(x){+}v^\top(w{-}x),\forall w{\in}\Re^n\right\},
\end{align}
which, in our case, is given by: $v_i=\lambda\sign(x_i)$, for $x_i\neq 0$ and $|v_i|\leq \lambda$, for $x_i=0$.
Therefore, the optimality conditions~\eqref{eq:gen_optCond} for the LASSO problem~\eqref{eq:str_prb} become
\begin{subequations}\label{eq:opt_LASSO}
\begin{align}
-\nabla_if(x^\star)=\lambda\sign(x^\star_i),\qquad \textrm{if}\  x^\star_i\neq 0,\label{eq:opt_LASSOa}\\
|\nabla_if(x^\star)|\leq\lambda,\qquad \textrm{if}\ x^\star_i= 0,\label{eq:opt_LASSOb}
\end{align}
\end{subequations}
Suppose for a moment that we knew the partition of indices $\alpha$ and $\beta$ corresponding to the nonzero and zero components of an optimal solution $x^\star$ respectively, as well as the signs of the nonzero components. Then we would be able to compute the nonzero components of $x^\star$ by solving the following linear system corresponding to~\eqref{eq:opt_LASSOa}:
\begin{equation}\label{eq:linsys_opt}
A_\alpha^\top A_\alpha x^\star_{\alpha}=A_\alpha^\top y-\lambda\sign(x^\star_\alpha).
\end{equation}
Notice that the support set $\alpha$ of $x^\star$ is much smaller compared to its dimension $n$. Hence, provided that $\alpha$ (as well as the signs of non-zero entries) has been identified, the problem becomes very easy. 
Roughly speaking, the algorithm we will develop can be interpreted as a fast procedure for automatically identifying the partition $\{\alpha,\beta\}$ corresponding to an optimal solution by solving a sequence of linear systems of the form~\eqref{eq:linsys_opt}. On the other hand, it can be seen as a Newton method for solving the following reformulation of the optimality conditions~\eqref{eq:gen_optCond}:  
\begin{equation}\label{eq:proxOpt}
x=\prox_{\gamma g}(x-\gamma\nabla f(x)),
\end{equation}
where $\prox_{\gamma g}$ is the \emph{proximal mapping} defined as
\begin{align}
 \prox_{\gamma g}(z)=\argmin_{x\in\Re^n}\left\{g(x)+\tfrac{1}{2\gamma}\|x-z\|_2^2\right\},
\end{align}
and $\gamma$ is taken smaller than the Lipschitz constant of $f$, i.e., $\gamma<1/\|A\|^2$. 
In the case where $g=\lambda\|\cdot\|_1$, $\prox_{\gamma g}$ is the soft-thresholding operator:
\begin{equation}\label{eq:SoftThres}
(\prox_{\gamma g}(z))_i=\sign(z_i)(|z_i|-\gamma\lambda)_+,\quad i=1,\ldots,n.
\end{equation}
The iterative soft-thresholding algorithm (ISTA) is a fixed point iteration for solving~\eqref{eq:proxOpt}. In fact, it is just an application of the well known forward-backward splitting technique for solving~\eqref{eq:gen_optCond}. On the other hand, FISTA is an accelerated version of ISTA where an extrapolation step between the current and the previous step precedes the forward-backward step~\cite{Beck2009,nesterov2013gradient}. 

To motivate our algorithm, instead of a fixed point problem, we view~\eqref{eq:proxOpt} as a problem of finding a zero of the so-called \emph{fixed point residual}:
\begin{equation}\label{eq:FPR}
R_\gamma(x):=x-\prox_{\gamma g}(x-\gamma\nabla f(x)).
\end{equation}
One then would be tempted to apply Newton's method for finding a root of~\eqref{eq:FPR}. Unfortunately the fixed point residual is not everywhere differentiable, hence the classical Newton method is not well-defined. However,
it is well known that $R_\gamma$ is nonexpansive, hence globally Lipschitz continuous~\cite{patrinosFBN}. Therefore, the machinery of nonsmooth analysis can be employed to devise a generalized Newton method for $R_\gamma(x)=0$, namely the \emph{semismooth} Newton method. Due to a celebrated theorem by Rademacher, Lipschitz continuity of $R_\gamma$ implies almost everywhere differentiability. Let $\mathcal{F}$ stand for the set of points where $R_\gamma$ is differentiable.
The $B$-differential  of the nonsmooth mapping $R_\gamma$ at $x$ is defined by
{
\small
\begin{align}
 \partial_B R_\gamma(x):=\left\{B\in\Re^{n\times n}\left| 
 \begin{array}{l}
 \exists\{x_n\}\in\mathcal{F}: x_n\to x,\\
 R'_\gamma(x_n)\to B
 \end{array}
  \right. \right\}.
\end{align}
}

If $R_\gamma$ is continuously differentiable at a point $x\in\Re^n$ then $\partial_B R_\gamma(x)=\{R'_\gamma(x)\}$. Otherwise $\partial_B R_\gamma$ may contain more than one elements (matrices in $\Re^{n\times n}$). 
The semismooth Newton method for solving $R_\gamma(x)=0$ is simply
\begin{equation}\label{eq:SNewton}
x^{k+1}=x^k-H_k^{-1}R_\gamma(x^k),\qquad H_k\in\partial_B R_\gamma(x^k).
\end{equation}
Since in the case of LASSO the fixed point residual is piecewise affine, it is also strongly semismooth. Provided that solution $x^\star$ of~\eqref{eq:str_prb} is unique (which is the case if for example the entries of $A$ are drawn Gaussian i.i.d~\cite{LASSO_unique}) and that the initial iterate $x^0$ is close enough to $x^\star$, the sequence of iterates defined by~\eqref{eq:SNewton} is well defined (any matrix in $\partial_B R_\gamma(x^k)$ is nonsingular) and converges to $x^\star$ at a quadratic rate, \emph{i.e.}, $\limsup_{k\to\infty}\frac{\|x^{k+1}-x^\star\|}{\|x^k-x^\star\|^2}<\infty$~\cite{patrinosFBN}. In fact, since the fixed point residual is piecewise affine for the LASSO problem, it can be shown that~\eqref{eq:SNewton} converges in a \emph{finite number of iterations}, in exact arithmetic.
Specializing iteration~\eqref{eq:SNewton} to the LASSO problem, an element $H_k$ of $\partial_B R_\gamma(x_k)$ takes the form
$H_k=I-P_k(I-\gamma A^\top A)$. Here $P_k$ is a diagonal matrix with $(P_k)_{ii}=1$ for $i\in\alpha_k$ and $(P_k)_{ii}=0$, for $i\in\beta_k$, where 
\begin{subequations}\label{eq:ab_k}
\begin{align}
\alpha_k&=\{i\mid |x^k_i-\gamma\nabla_i f(x^k)|>\gamma\lambda\},\label{eq:alphak}\\
\beta_k&=\{i\mid |x^k_i-\gamma\nabla_i f(x^k)|\leq \gamma\lambda\}.
\end{align}
\end{subequations}
Computing the Newton direction amounts to solving the so-called \emph{Newton system}, $H_kd^k=-R_\gamma(x^k)$. 
Taking advantage of the special structure of $P_k$ and applying some permutations, this simplifies to 
\begin{subequations}\label{eq:NewtonDir}
\begin{align}
d^k_{\beta_k}&=-(R_\gamma(x^k))_{\beta_k},\\
\gamma A_{\alpha_k}^\top A_{\alpha_k}d^k_{\alpha_k}&=-(R_\gamma(x^k))_{\alpha_k}-\gamma A_{\alpha_k}^\top A_{\beta_k}d_{\beta_k}^k.
\end{align}
\end{subequations}
Taking into consideration~\eqref{eq:SoftThres}, we obtain
\begin{align}
 (R_\gamma(x))_i=\begin{cases}\gamma(\nabla_if(x)+s_i(x)\lambda),& i\in\alpha,\\ x_i,& i\in\beta,\end{cases}
\end{align}
where $s_i(x)=\sign(x_i-\gamma\nabla_if(x))$, and $\alpha,\beta$ are defined using~\eqref{eq:NewtonDir} (dropping index $k$). Therefore, after further rearrangement and simplification, the Newton system becomes
\begin{subequations}\label{eq:NewtonDirSimple}
\begin{align}
d^k_{\beta_k}&=-x^k_{\beta_k},\\
 A_{\alpha_k}^\top A_{\alpha_k}(x^k_\alpha+d^k_{\alpha_k})&=A_{\alpha_k}^\top y -\lambda s_{\alpha_k}(x^k)\label{eq:linsysNewt}.
\end{align}
\end{subequations}
Summing up, an iteration of the semismooth Newton method for solving the LASSO problem takes the form $x^{k+1}=x^k+d^k$, which becomes: 
\begin{subequations}\label{eq:LassoNewton}
\begin{align}
x^{k+1}_{\beta_k}&=0,\\
x^{k+1}_{\alpha_k}&=( A_{\alpha_k}^\top A_{\alpha_k})^{-1}(A_{\alpha_k}^\top y -\lambda s_{\alpha_k}(x^k)).
\end{align}
\end{subequations}

Note a connection with~\eqref{eq:linsys_opt}: the index sets $\alpha_k$, $\beta_k$ serve as estimates for the nonzero and zero components of $x^\star$. 

There are two obstacles that we need to overcome before we arrive at a sound, globally convergent algorithmic scheme for solving the LASSO problem. The first one is that the semismooth Newton method~\eqref{eq:LassoNewton} converges only when started close to the solution. 
The second obstacle concerns the fact that~\eqref{eq:LassoNewton} might not be well defined, in the sense that  $A_{\alpha_k}$ might not have full column rank, hence  $A_{\alpha_k}^\top A_{\alpha_k}$ will be singular. Indeed when $\lambda$ is very small or when $x^k$ is far from the solution, the index set $\alpha_k$ defined in~\eqref{eq:alphak} might have cardinality larger than $m$ (which is the only case where a singularity may arise given our construction of $A$), the number of rows of $A$.  The next two subsections are devoted to proposing strategies to overcome these two issues.

Overall, the algorithm can be seen as an active set strategy where large changes 
on the active set are allowed in every iteration (instead of only one index) 
leading to faster convergence.

\subsection{Globalization strategy}
To enforce global convergence of Newton-type methods for solving nonlinear systems of equations it is customary to use a merit function based on which a step $\tau_k$ is selected which guarantees that 
\begin{equation}\label{eq:DampedNewton}
x^{k+1}=x^k+\tau_kd^k
\end{equation}
decreases the merit function sufficiently. 

Recently the following merit function, namely the \textit{Forward-Backward Envelope} (FBE), 
was proposed for problems of the form~\eqref{eq:str_prb}
\small
\begin{equation}\label{eq:FBE}
\varphi_\gamma(x)=\inf_{z\in\Re^n}\{f(x)+\nabla f(x)^\top(z-x)+g(z)+\tfrac{1}{2\gamma}\|z-x\|_2^2\}.
\end{equation}
\normalsize
It is easy to see that function inside the infimum is strongly convex with respect to $z$ and the infimum is uniquely achieved by the forward-backward step $T_\gamma(x)=\prox_{\gamma g}(x-\gamma\nabla f(x))$. Therefore, in order to evaluate $\varphi_\gamma$ at a point $x$ one simply needs to be be able to perform the same operations required by FISTA. Furthermore, function $\varphi_\gamma$ is continuously differentiable with gradient given by
\begin{align}
 \nabla\varphi_\gamma(x)=\gamma^{-1}(I-\gamma\nabla^2f(x))R_\gamma(x).
\end{align}
If $\gamma<1/\|A\|^2$, then solutions of~\eqref{eq:str_prb} are exactly the stationary points, \emph{i.e.} the points for which $\nabla\varphi_\gamma(x)$ becomes zero. In fact one can additionally show that  minimizing the FBE (which is an unconstrained smooth optimization problem) is entirely equivalent to solving~\eqref{eq:str_prb}, in the sense that
$\inf\ \varphi=\inf \varphi_\gamma$ and $\argmin\ \varphi=\argmin\ \varphi_\gamma$.
In the case of LASSO where $f$ is quadratic the FBE is convex.

It is not hard to check that, provided it is well defined, the semismooth Newton direction $d^k$ given by~\eqref{eq:NewtonDir} is in fact a direction of descent for the FBE, \emph{i.e.}, $\nabla\varphi_\gamma(x^k)^\top d^k<0$. Therefore, one can perform a standard backtracking line-search to find a suitable step that guarantees the Armijo condition and hence global convergence: Pick the first nonnegative integer $i_k$ such that $\tau_k=2^{-i_k}$ satisfies
\begin{equation}\label{eq:Armijo}
\varphi_\gamma(x_k+\tau_kd_k)\leq\varphi_\gamma(x_k)+\zeta\tau_k\nabla\varphi_\gamma(x_k)^\top d_k.
\end{equation}
where  $\zeta\in(0,1/2)$, and then set $x_{k+1}=x_k+\tau_kd_k$. Furthermore, this step choice guarantees that as soon as $x^k$ is close enough to the solution, $\tau_k=1$ will always satisfy~\eqref{eq:Armijo} and the iterates will be given by the (pure) semismooth Newton method~\eqref{eq:LassoNewton} inheriting all its convergence properties.

\subsection{Continuation strategy}
Since the goal of solving the LASSO problem is to recover the sparsest solution of 
$Ax=y$, we know that such a solution will have a $\|x^\star\|_0$ much smaller than $m$, 
hence as soon as $x^k$ is  close to the solution, $A_{\alpha_k}$ will have full column rank.

When $\lambda$ is small or when $x^0$ is far from the solution, the cardinality 
of~\eqref{eq:alphak} can be larger than $m$ rendering the left-hand side matrix in~\eqref{eq:linsysNewt} singular. 
This said, we know that iterates for which the set $\alpha_k$ contains more than $m$
indices, are known to be away from the optimal solution since we know for sure that at 
the solution the cardinality of
$\alpha(x^\star)=\{i\mid |x^\star_i-\gamma\nabla f(x^\star)|>\gamma\lambda\}=\{i\mid x^\star_{i}\neq 0\}$ 
will be (much) smaller than $m$. 

{In order to enforce that $\alpha_k$ contains few elements, a simple continuation strategy 
that gradually reduces $\lambda_k$ to the target value $\lambda$ is employed.
Following~\cite{xiao2013proximal}, we start with $\lambda_0=\max\{\lambda, \|A^{\top}\nabla f(x_0)\|_{\infty}\}$
and we decrease $\lambda_k = \max\{\lambda, \eta \lambda_{k-1}\}$, for some $\eta\in(0,1)$ whenever 
$\|R_\gamma(x_k)\| \leq \lambda_k \epsilon_k$ with $\epsilon_k\to 0$.
Not only does this ensure that~\eqref{eq:NewtonDir} is well defined, but it allows to solve linear 
systems of small dimension to determine the Newton direction. 
A conceptual pseudo-algorithm summarizing the basic steps of the proposed algorithm
is shown hereafter.}

\begin{algorithm}[!ht]
	\caption{\textbf{Forward-Backward Newton with continuation}}
	\begin{algorithmic}
		\Require $A$, $y$, $x^0\in\Re^n$ (initial guess), $\gamma\in(0,1/\|A\|^2)$, $\lambda > 0$,  $\eta\in(0,1)$, $\epsilon$ (tolerance)	
		\State $\bar{\lambda} \leftarrow \max\{\lambda, \|A^{\top}\nabla f(x_0)\|_{\infty}\}$, $\bar{\epsilon} \leftarrow \epsilon$
		\While{$\bar{\lambda}  > \lambda$ or $\|R_{\gamma}(x^k;\bar{\lambda})\|>\epsilon$}		  
		  \State $x^{k+1}=x^k+\tau_kd^k$, 
		  \State where $d^k$ solves~\eqref{eq:NewtonDirSimple} and $\tau_k$ satisfies~\eqref{eq:Armijo}
		  \If{$\|R_\gamma(x^k;\bar{\lambda})\|  \leq  \lambda \bar{\epsilon}$}
		    \State $\bar{\lambda} \leftarrow \max\{ \lambda, \eta \bar{\lambda}\}$, $\bar{\epsilon}\leftarrow \eta \bar{\epsilon}$
		  \EndIf
		\EndWhile
	\end{algorithmic}
	\label{alg:recursivecs}
\end{algorithm}

In Algorithm~\ref{alg:recursivecs} we denote by $R_\gamma(x; \bar\lambda)$ the 
fixed point residual introduced in~\eqref{eq:FPR} replacing $g$ by $\bar{g}(x) = \bar{\lambda}\|x\|_1$.

\section{Simulations}\label{sec:Simulations}
In this section we apply the proposed methodology to various data streams and 
we compare it to standard algorithms such as FISTA, ADMM~\cite{boyd2011distributed} and the interior point
method of Kim \textit{et al.}~\cite{kim2007interior}, aka L1LS method. 
In our approach we used the continuation strategy described above with $\eta=0.5$
and an Armijo line search.
The required tolerance for the termination of all algorithms was set to $10^{-8}$.

We observed that after decompressing the first window, the number of iterations 
required for convergence was remarkably low (in most cases, around $4$ iterations 
for each window were sufficient). It should be highlighted that after the first decompression, 
the computational cost of the algorithm decreases 
significantly. This is first because of the aforementioned warm-start and second 
because the value of the residual $r = Ax-y$ is updated by using only vector-vector operations.
Updating $A^{(i+1)}\leftarrow A^{(i)}P$ and $x_0^{(i+1)} \leftarrow P^{\top}\hat{x}^{(i)}$ and using 
the fact that $P$ is orthogonal, we have that
$r^{(i+1)} = A^{(i+1)}x_0^{(i+1)}-y^{(i+1)} = r^{(i)} + y^{(i)} - y^{(i+1)}$.

A stream of total length $N=10^6$ was generated as follows: its entries are drawn from  $\mathrm{Ber}(S)$, where $S$ is average stream sparsity. 
Then, the non-zero entries are taken uniform in $[-2,-1] \cup [1,2]$ and multiplied by $8 \sigma \sqrt{2 \log N}$ 
(dynamic range assumption~\cite{candes2009}), based on selected noise variance $\sigma^2$.  
We then select window size $n$, let $s=nS, m=4s$, and generate sampling matrix  $A$ with i.i.d. entries $\mathcal{N}(0, 1/m)$. 
For LASSO, we pick $\lambda = 4 \sigma \sqrt{2 \log n}$ following~\cite{candes2009}.

\begin{figure}
 \centering
 \includegraphics[keepaspectratio=true,width=0.46\textwidth]{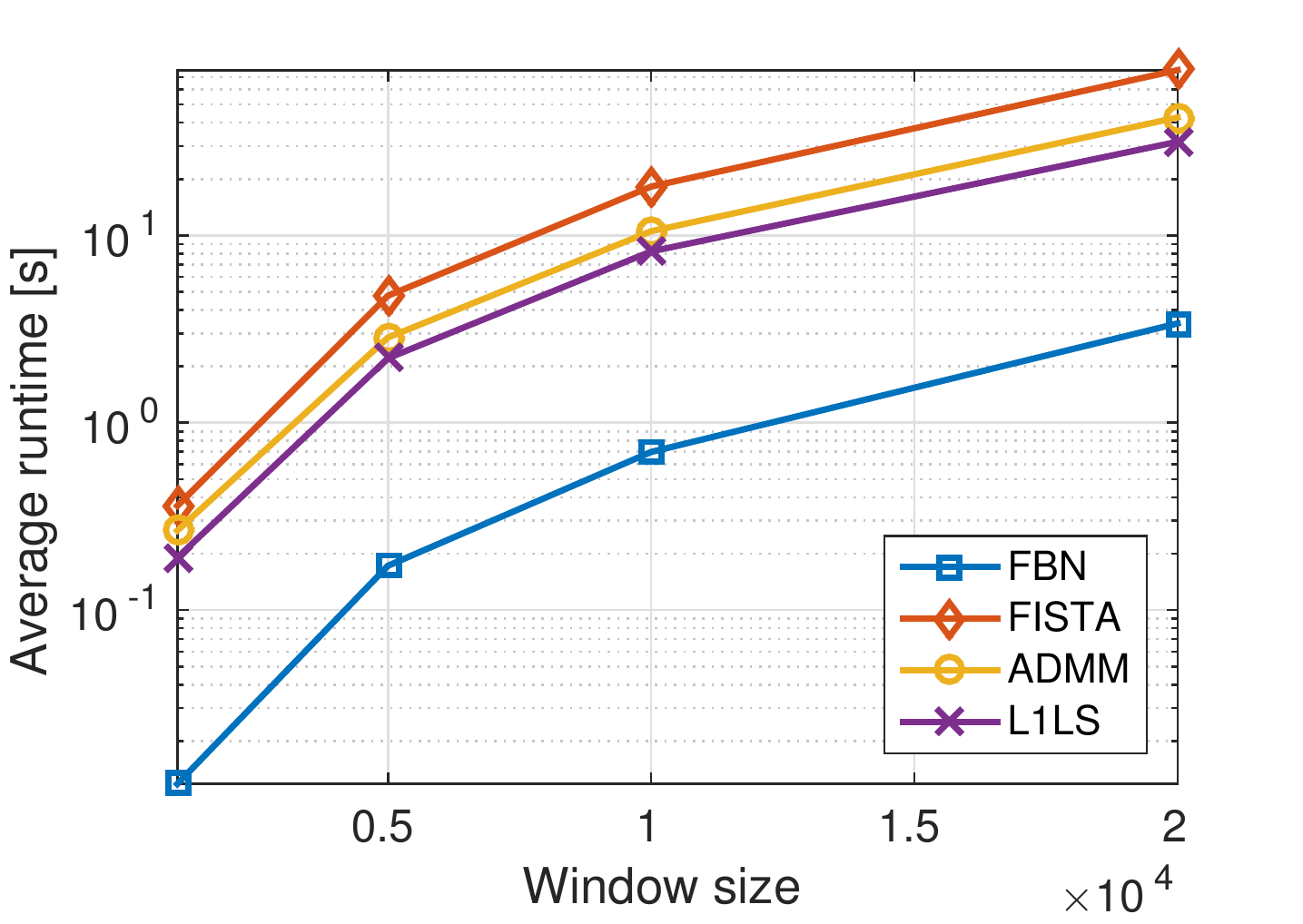}
 \caption{Average runtimes varying the window size. FBN: proposed forward-backward 
          Newton method, L1LS~\cite{kim2007interior}, ADMM~\cite{Afonso2010}, FISTA~\cite{Beck2009}.}
 \label{fig:runtimes_n}
\end{figure}

\begin{figure}
 \centering
 \includegraphics[keepaspectratio=true,width=0.46\textwidth]{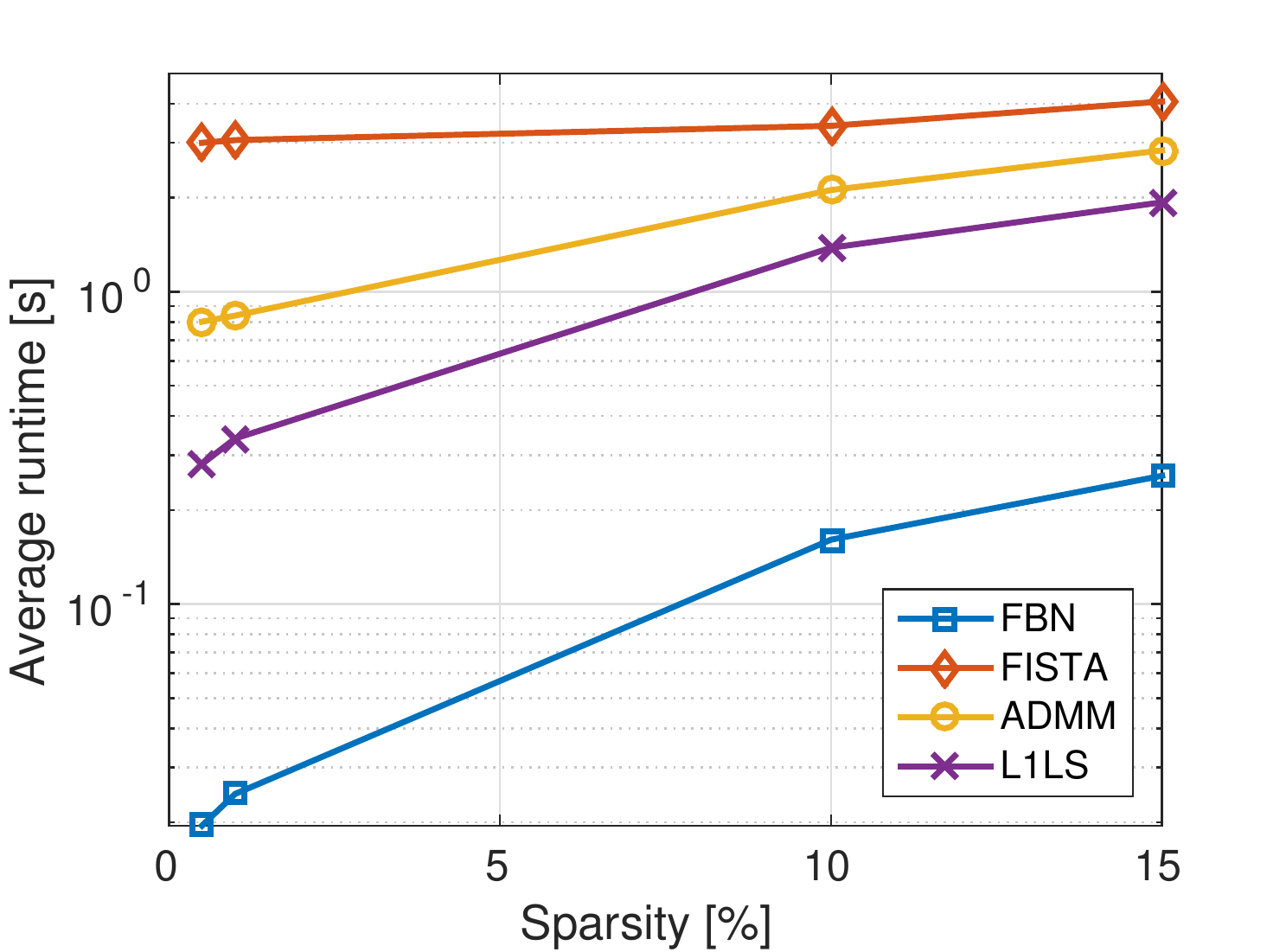}
 \caption{Average runtimes varying the sparsity of the data stream.}
 \label{fig:runtimes_sparsity}
\end{figure}

In Figures~\ref{fig:runtimes_n} and~\ref{fig:runtimes_sparsity} 
we observe that the proposed algorithm outperforms all state-of-the-art 
methods by an order of magnitude.
The results presented in Figure 1 were obtained for a fixed sparsity $10\%$ and
$w^{(i)}$ being a zero-mean normally distributed noise with variance $0.01$.
In Figure 2, for the same noise level and using a window of size $n=5000$ 
we show how the runtime is affected by the sparsity of the data 
stream\footnote{We have also conducted extensive experiments varying SNR, termination accuracy 
and comparing the speed of RCS against classical CS; in all of them we have observed significant improvements, 
but we exclude them for length considerations.}.

\section{Conclusions}\label{sec:Conclusions}

We have proposed an efficient method for successively decompressing the entries of a data stream sampled using Recursive Compressed Sensing~\cite{RCS}.
We have proposed a second-order proximal method for solving LASSO with accelerated convergence over state-of-art methods. 
Our scheme is very efficient as each iteration nails down to solving a linear system of low dimension, which we may further avoid by a single Cholesky factorization at a pre-processing step. We have tested our algorithm against the state-of-art for various windows sizes and sparsity patterns; our experiments depict notable speed-up which may renders RCS suitable for an online implementation under stringent time constraints.


\begin{spacing}{0.95}
\bibliographystyle{IEEEtran}
\bibliography{EUSIPCObib}

\begin{thebibliography}{10}
\providecommand{\url}[1]{#1}
\csname url@samestyle\endcsname
\providecommand{\newblock}{\relax}
\providecommand{\bibinfo}[2]{#2}
\providecommand{\BIBentrySTDinterwordspacing}{\spaceskip=0pt\relax}
\providecommand{\BIBentryALTinterwordstretchfactor}{4}
\providecommand{\BIBentryALTinterwordspacing}{\spaceskip=\fontdimen2\font plus
\BIBentryALTinterwordstretchfactor\fontdimen3\font minus
  \fontdimen4\font\relax}
\providecommand{\BIBforeignlanguage}[2]{{%
\expandafter\ifx\csname l@#1\endcsname\relax
\typeout{** WARNING: IEEEtran.bst: No hyphenation pattern has been}%
\typeout{** loaded for the language `#1'. Using the pattern for}%
\typeout{** the default language instead.}%
\else
\language=\csname l@#1\endcsname
\fi
#2}}
\providecommand{\BIBdecl}{\relax}
\BIBdecl

\bibitem{RCS}
\BIBentryALTinterwordspacing
N.~Freris, O.~\"{O}\c{c}al, and M.~Vetterli, ``{Recursive Compressed
  Sensing},'' Tech. Rep., 2014. [Online]. Available:
  \url{http://arxiv.org/abs/1312.4895}
\BIBentrySTDinterwordspacing

\bibitem{Donoho2006}
D.~L. Donoho, ``{Compressed sensing},'' \emph{IEEE Transactions on Information
  Theory}, vol.~52, pp. 1289--1306, 2006.

\bibitem{Candes2006b}
E.~Cand{\`e}s and T.~Tao, ``Near-optimal signal recovery from random
  projections: Universal encoding strategies?'' \emph{IEEE Transactions on
  Information Theory}, vol.~52, no.~12, pp. 5406--5425, 2006.

\bibitem{Candes2008}
E.~Cand\`{e}s and M.~Wakin, ``An introduction to compressive sampling,''
  \emph{IEEE Signal Processing Magazine}, vol.~25, no.~2, pp. 21--30, 2008.

\bibitem{CS_applications}
S.~Qaisar, R.~Bilal, W.~Iqbal, M.~Naureen, and S.~Lee, ``Compressive sensing:
  From theory to applications, a survey,'' \emph{Journal of Communications and
  Networks}, vol.~15, no.~5, pp. 443--456, 2013.

\bibitem{Tibshirani1996}
R.~Tibshirani, ``{Regression Shrinkage and Selection via the LASSO},''
  \emph{Journal of the Royal Statistical Society. Series B (Methodological)},
  vol.~58, pp. 267--288, 1996.

\bibitem{candes2009}
E.~Cand{\`e}s and Y.~Plan, ``Near-ideal model selection by $\ell_1$
  minimization,'' \emph{The Annals of Statistics}, vol.~37, pp. 2145--2177,
  2009.

\bibitem{kim2007interior}
S.-J. Kim, K.~Koh, M.~Lustig, S.~Boyd, and D.~Gorinevsky, ``An interior-point
  method for large-scale {$\ell_1$}-regularized least squares,'' \emph{IEEE
  Journal of Selected Topics in Signal Processing}, vol.~1, no.~4, pp.
  606--617, 2007.

\bibitem{Parikh:2013}
N.~Parikh and S.~Boyd, ``{Proximal algorithms},'' \emph{Foundations and Trends
  in Optimization}, vol.~1, no.~3, pp. 123--231, 2013.

\bibitem{Beck2009}
A.~Beck and M.~Teboulle, ``A fast iterative shrinkage-thresholding algorithm
  for linear inverse problems,'' \emph{SIAM Journal on Imaging Sciences},
  vol.~2, no.~1, pp. 183--202, 2009.

\bibitem{Wright:2009kw}
S.~J. Wright, R.~D. Nowak, and M.~A.~T. Figueiredo, ``{Sparse Reconstruction by
  Separable Approximation},'' \emph{IEEE Transactions on Signal Processing},
  vol.~57, no.~7, pp. 2479--2493, 2009.

\bibitem{Afonso2010}
M.~Afonso, J.~Bioucas-Dias, and M.~A.~T. Figueiredo, ``Fast image recovery
  using variable splitting and constrained optimization,'' \emph{IEEE
  Transactions on Image Processing}, vol.~19, no.~9, pp. 2345--2356, 2010.

\bibitem{allerton}
N.~Freris, O.~\"{O}\c{c}al, and M.~Vetterli, ``{Compressed Sensing of Streaming
  data},'' in \emph{Proceedings of the 51st Allerton Conference on
  Communication, Control and Computing}, pp. 1242--1249.

\bibitem{patrinos2013proximal}
P.~Patrinos and A.~Bemporad, ``Proximal {N}ewton methods for convex composite
  optimization,'' in \emph{IEEE Conference on Decision and Control}, Florence,
  Italy, 2013, pp. 2358--2363.

\bibitem{patrinosFBN}
\BIBentryALTinterwordspacing
P.~Patrinos, L.~Stella, and A.~Bemporad, ``Forward-backward truncated {N}ewton
  methods for convex composite optimization,'' Tech. Rep., 2014. [Online].
  Available: \url{http://arxiv.org/abs/1402.6655}
\BIBentrySTDinterwordspacing

\bibitem{nesterov2013gradient}
Y.~Nesterov, ``Gradient methods for minimizing composite functions,''
  \emph{Mathematical Programming}, vol. 140, no.~1, pp. 125--161, 2013.

\bibitem{LASSO_unique}
\BIBentryALTinterwordspacing
R.~J. Tibshirani, ``The {L}asso problem and uniqueness,'' Tech. Rep., 2012.
  [Online]. Available: \url{http://arxiv.org/abs/1206.0313}
\BIBentrySTDinterwordspacing

\bibitem{xiao2013proximal}
L.~Xiao and T.~Zhang, ``A proximal-gradient homotopy method for the sparse
  least-squares problem,'' \emph{{SIAM} Journal on Optimization}, vol.~23,
  no.~2, pp. 1062--1091, 2013.

\bibitem{boyd2011distributed}
S.~Boyd, N.~Parikh, E.~Chu, B.~Peleato, and J.~Eckstein, ``Distributed
  optimization and statistical learning via the alternating direction method of
  multipliers,'' \emph{Foundations and Trends in Machine Learning}, vol.~3,
  no.~1, pp. 1--122, 2011.

\end{thebibliography}
\end{spacing}

\end{document}